\documentclass[12pt]{amsart}

\usepackage[margin=1in, includefoot, footskip=30pt]{geometry}
\usepackage{amsmath, amsfonts, algorithm, booktabs, enumerate}
\usepackage[noend]{algpseudocode}
\usepackage[bookmarks=true]{hyperref}
\usepackage{setspace, wrapfig, graphicx}
\DeclareMathOperator*{\argmin}{arg\,min}
\DeclareMathOperator*{\argmax}{arg\,max}
\newcommand{\smfrac}[2]{\mbox{$\frac{#1}{#2}$}}

\algnewcommand\algorithmicinput{\textbf{Input:}}
\algnewcommand\Input{\item[\algorithmicinput]}
\algnewcommand\algorithmicoutput{\textbf{Output:}}
\algnewcommand\Output{\item[\algorithmicoutput]}
\newtheorem{theorem}{Theorem}[section]
\newtheorem{definition}[theorem]{Definition}

\title{Solving constrained quadratic binary problems\\via quantum adiabatic
evolution}
\date{\today}


\author{Pooya Ronagh}
\author{Brad Woods}
\author{Ehsan Iranmanesh}

\address[Pooya Ronagh, Brad Woods, and Ehsan Iranmanesh]{
1QBit,
458--550 Burrard Street,
Vancouver, British Columbia,
Canada \, V6C 2B5}
\email[Pooya Ronagh]{pooya.ronagh@1qbit.com}
\email[Brad Woods]{brad.woods@1qbit.com}
\email[Ehsan Iranmanesh]{ehsan.iranmanesh@1qbit.com}

\begin{document}

\begin{abstract}
Quantum adiabatic evolution is perceived as useful for binary quadratic
programming problems that are a priori unconstrained. For constrained problems,
it is a common practice to relax linear equality constraints as penalty terms in
the objective function. However, there has not yet been proposed a method for
efficiently dealing with inequality constraints using the quantum adiabatic
approach. In this paper, we give a method for solving the Lagrangian dual of a
binary quadratic programming (BQP) problem in the presence of inequality
constraints and employ this procedure within a branch-and-bound framework for
constrained BQP (CBQP) problems.
\end{abstract}

\maketitle
\onehalfspacing
\let\thefootnote\relax\footnotetext{
\hspace*{2.7mm} \emph{Published in:} Quantum Information \& Computation,
{Vol.~16}, {No.~11\&12}, {1029--1047}, Rinton~Press (2016).}

\section{Introduction}
An unconstrained binary quadratic programming (UBQP) problem is defined by
\begin{eqnarray} 
 &\mbox{Minimize} & x^TQx \nonumber \\
 & \text{subject to} & x \in \{0, 1\}^n, 
 \end{eqnarray} 
where, without loss of generality, $Q\in \mathbb{Z}^{n\times n}$. Recent
advancements in quantum computing technology~\cite{Boxio2015, DwaveNature,
GoogleNature} have raised hopes of the production of computing systems that are
capable of solving UBQP problems, and showing quantum speedup. The stochastic
nature of such systems, together with extant sources of noise and error, are
challenges yet to be overcome in achieving scalable quantum computing systems of
this type. This paper is nevertheless motivated by the assumption of the
existence of systems that can solve UBQP problems efficiently and to optimality,
or at least in conjunction with a framework of noise analysis of the suboptimal
results. We call such a computing system a UBQP oracle.

Many NP-hard combinatorial optimization problems arise naturally or can easily
be reformulated as UBQP problems, such as the quadratic assignment problem, the
maximum cut problem, the maximum clique problem, the set packing problem, and
the graph colouring problem (see, for instance, Boros and
Pr{\'e}kopa~\cite{Boros:1989}, Boros and Hammer~\cite{Boros:2002}, Bourjolly et
al.~\cite{Bourjolly:1994}, Du and Pardalos~\cite{Du:1999}, Pardalos and
Rodgers~\cite{Pardalos:1990,Pardalos:1992}, Pardalos and
Xue~\cite{Pardalos:1994}, and Kochenberger et al.~\cite{Kochenberger:2004}).

Numerous interesting applications that are expressed naturally in the form of
UBQP problems appear in the literature. Barahona et al.~\cite{Barahona:1988,
deSimone:1995} formulate and solve the problem of finding exact ground states of
spin glasses with magnetic fields. Alidaee et al.~\cite{Alidaee:1994} study the
problem of scheduling $n$ jobs non-preemptively on two parallel identical
processors to minimize weighted mean flow time as a UBQP problem. Bomze et
al.~\cite{Bomze:1999} give a comprehensive discussion of the maximum clique (MC)
problem. Included is the UBQP representation of the MC problem and a variety of
applications from different domains. UBQP has been used in the prediction of
epileptic seizures~\cite{Iasemidis:2001}. Alidaee et al.~\cite{Alidaee:2005}
discuss a number partitioning problem, formulating a special case as a UBQP
problem.

In this paper, we consider constrained binary quadratic programming (CBQP)
problems with linear constraints, stated formally as
\begin{eqnarray}
	\mathrm{(P)} & \mbox{Minimize} & x^TQx \nonumber\\
	& \mbox{subject to} & Ax\leq b, \nonumber\\
	& & x\in \{ 0,1\} ^n\,,
\end{eqnarray}
where $Q\in \mathbb{Z}^{n\times n}$ and $A\in \mathbb{Z}^{m\times n}$.

There are many problems which naturally occur as linearly constrained binary
quadratic programming problems. To illustrate, consider the well-studied
quadratic assignment problem: the problem of assigning facilities to locations
where the cost is a function of the distance and flow between facilities plus
the cost of assigning a facility to a specific location. The problem requires
that each facility be assigned exactly one location, and each location exactly
one facility, and is easily expressed in the form of
(P)~\cite{lawler1963quadratic}. Other examples include the clique partitioning
problem~\cite{oosten2001clique}, the quadratic minimum spanning tree
problem~\cite{assad1992quadratic}, and the quadratic shortest path
problem~\cite{rostami2015quadratic}.

In the literature, CBQP problems are commonly reformulated as UBQP problems by
including quadratic penalties in the objective function as an alternative to
explicitly imposing constraints. Although this method has been used very
successfully on classical hardware, it is not a viable approach when using
quantum adiabatic hardware, as the reformulation dramatically increases the
density, range of coefficients, and dimension of the problem.

We present a branch-and-bound approach which uses Lagrangian duality to solve
(P) and show that a UBQP oracle can be used to solve the Lagrangian dual (LD)
problem with successive applications of linear programming (LP). Throughout
this paper, we will refer to this algorithm as the quantum branch-and-bound
algorithm. We introduce the notion of \emph{quantum annealing leniency}, which
can be used to compare a classical algorithm running on a Turing machine to an
algorithm running on an oracle Turing machine \cite{Sipser} with a quantum
annealing oracle. This measure represents the maximum threshold of the average
time an oracle query is allowed to take in order to outperform the benchmark
algorithm running on a classical Turing machine. In our experiments, we
benchmark our quantum branch-and-bound algorithm against the Gurobi Optimizer
\cite{gurobi}. The quantum annealing leniency of the quantum branch-and-bound
algorithm with respect to the Gurobi Optimizer is measured.

This paper is organized as follows. Section \ref{sec:Prelim} gives an overview
of the quantum adiabatic approach to computation and the relationship of the
approach to UBQP. Section \ref{sec:LB} presents lower- bounding procedures for
CBQP. Section \ref{sec:UB} presents a local search heuristic for CBQP.
Branching strategies are described in Section \ref{sec:Branching}. All of these
algorithms are then integrated in the quantum branch-and-bound framework
presented in Section \ref{sec:qbnb}. Test instances and results of our
computational experiments are presented in Section~\ref{sec:Experiments}. In
Section~\ref{sec:Discussion}, we provide practical instruction for the
application of our method to the quantum processors manufactured by D-Wave
Systems Inc., and in Section~\ref{sec:Generalization}, we mention how our method
can more generally be used to solve constrained binary programming problems with
higher-order polynomial objectives and higher-order polynomial constraints.

\section{Computing using quantum adiabatic evolution}\label{sec:Prelim}

Recent advancements in quantum hardware technology have motivated an increase in
the study of forms of computation that differ in computational complexity from
Turing machines. The quantum gate model is a means of achieving powerful quantum
algorithms such as Shor's well known quantum algorithms for integer
factorization and computing discrete logarithms~\cite{shor}. Aside from the
quantum gate model, there are several other paradigms of quantum information
technology, each of which would open a new world of possible algorithm designs
to be realized on a corresponding practical quantum processor.

Farhi et al. \cite{farhi00,farhi01} propose quantum adiabatic evolution as a
novel paradigm for the design of quantum algorithms. Quantum adiabatic
computation is expressed by the Schr\"odinger equation of a time-dependent
Hamiltonian
\begin{eqnarray}\label{eq:adiabatic-ev}
H(t) & := & \left(1-\smfrac{t}{T}\right)H_0
+ \smfrac{t}{T} H_f, \quad 0 \leq t \leq T\,.
\end{eqnarray}
Here $T$ is a constant \emph{delay factor}. The system is evolved according to
(\ref{eq:adiabatic-ev}) from an initial Hamiltonian $H_0$ at time $t=0$ to a
final Hamiltonian $H_f$ at time $t=T$. The former Hamiltonian is such that
setting the system to its ground state is easy, and the latter Hamiltonian is
constructed from a polynomial objective function $f(z)$ in binary variables.
$H_f$ is associated to $f$ such that the range of $f$ is identical to the
eigenvalue spectrum of $H_f$. By the quantum adiabatic theorem~\cite{messiah},
when the system is initially set to the ground state of $H_0$ at time $t=0$, and
$T$ is sufficiently large, the system tends to stay in the ground state of
$H(t)$ for all $t$.

Van Dam et al.~\cite{vazirani} show that it is sufficient to have $T\in
O({\Delta_\mathrm{max}}/{g_\mathrm{min}^2}).$ Here $g_\mathrm{min}$ is the
minimum difference over time $s$ between the smallest two eigenvalues of $H(s)$
and $\Delta_\mathrm{max} := \max_s {\|{\smfrac{d}{ds}}H(s)\|_2}$. They give an
example of an adiabatic quantum algorithm for searching that matches the
quadratic speedup obtained by Grover's search algorithm. This example
demonstrates that the ``quantum local search,'' which is implicit in the
adiabatic evolution, is truly non-classical in nature from a computational
perspective. Also \cite[Theorem~1]{vazirani} explains how the continuous-time
evolution of $t \in [0, T]$ can be approximated by a quantum circuit consisting
of a sequence of $\mathrm{poly}(nT)$ unitary transformations.

All of the above considerations suggest that practical quantum hardware can
yield a significant quantum speedup in certain integer programming problems. Our
goal is to design and analyze optimization algorithms that work in conjunction
with such integer programming oracles. Specifically, we work under the
assumption of the existence of a \emph{UBQP oracle}, an oracle Turing machine
for solving UBQP problems. This assumption is motivated by prototypes of quantum
annealers recently manufactured by D-Wave Systems, where couplings connect pairs
of quantum bits \cite{DwaveNature}. Our suggested methods are easily
generalizable to take advantage of systems with higher-degree interactions of
quantum bits if such systems are implemented in the future (see
Section~\ref{sec:Generalization}).

As a final remark, it is important to mention that quantum annealers are coupled
to an environment, and this significantly affects their performance.
Albash~et~al.~\cite{albash2014} propose a noise model for D-Wave devices. This
model includes the control noise on the local field and couplings of the chip,
as well as the effect of the cross-talk between qubits that are not coupled.
Eventually \cite{albash2014} concludes that despite the thermal excitations and
small value of the ratio of the single-qubit decoherence time to the annealing
time, an open-system, quantum-dynamical description of the D-Wave device that
starts from a quantized energy level structure is well justified. The design of
benchmark instances that can detect quantum speedup or any quantum advantage of
a quantum annealer in comparison to state-of-the-art classical algorithms is
studied by Katzgraber~et~al.~\cite{katzgraber2015}. Zhu~et~al.~\cite{zoshk2016}
show that increasing the classical energy gap beyond the intrinsic noise level
of the machine can improve the success of the D-Wave Two quantum annealer, at
the cost of producing considerably easier benchmark instances. We refer the
reader to \cite{caty} for the practicality and best practices in using D-Wave
devices.

\section{Lower-bounding procedures}\label{sec:LB}

\subsection{Linearization relaxation}
A standard linearization of (P) involves relaxing the integrality constraint on
variables $x_i, (i= 1, \ldots, n)$, and defining continuous variables $y_{ij}$
for every pair $x_ix_j$ in the objective function with $i < j$, yielding the
following linearized problem.
\begin{eqnarray}
	\mathrm{(P_{LP})}& \mbox{  Minimize   } &
	\sum\limits_{1\leq i<j\leq n}  2\, q_{ij} y_{ij} + \sum_{i=1}^n  q_{ii} x_i
	\nonumber\\
		&\mbox{  subject to   } & Ax\leq b\,, \nonumber\\
	&& y_{ij} \geq x_i + x_j -1 \quad(\forall\, i,j
		\text{ such that }i<j \text{ and } q_{ij} >0)\,, \nonumber\\
	&& y_{ij} \leq x_i \quad(\forall\, i,j
		\text{ such that }i<j \text{ and } q_{ij} <0)\,,\nonumber \\
	&& y_{ij} \leq x_j \quad(\forall\, i,j
		\text{ such that }i<j \text{ and } q_{ij} <0)\,,\nonumber \\
	&& 0 \leq x_i \leq 1 \quad (i=1, \ldots, n)\,,\nonumber \\
	& & y\geq 0
\end{eqnarray}
A lower bound to (P) can now be obtained by solving $\mathrm{(P_{LP})}$ using
linear programming. We employed this linearization in our computational
experiments (see Section~\ref{sec:Experiments}). Note that there are several
methods for linearizing (P), many of which have been mentioned in the survey by
Floudas and Gounaris \cite{Floudas}; in this paper, however, we consider
$\mathrm{(P_{LP})}$ to be \emph{the} LP relaxation of (P).
 
\subsection{Lagrangian dual}
We can give a lower bound for (P) via the LD problem
\begin{eqnarray}
	\mathrm{(L)} & \max \limits_{\lambda \in \mathbb{R}^m_+} & d(\lambda) \,,
\end{eqnarray}
where $d(\lambda)$ is evaluated via the \textit{Lagrangian relaxation}
\begin{eqnarray}
\mathrm{(L_{\lambda})} & d(\lambda)
= \min \limits_{x\in \{ 0,1\}^n} L(x,\lambda)
= & x^TQx + \lambda^T(Ax-b)\,.
\end{eqnarray}
The function $d(\lambda)$ is the minimum of a finite set of linear functions of
$\lambda$ and hence it is concave and piecewise linear. The following theorem
shows that $\mathrm{L_{\lambda}}$ is a lower bound for (P). For any problem $Q$,
let $v(Q)$ be the optimal objective value.

\vspace*{12pt}
\noindent
{\bf Proposition~1 (Weak Duality):} For all $\lambda \in \mathbb{R}^m_+$, we
have $v(L_{\lambda}) \le v(P)\,.$

\vspace*{12pt}
\noindent
{\bf Proof:} Since every feasible solution for (P) is feasible for $L_\lambda$,
then for any $\lambda \in \mathbb{R}^m_+$ we have
\begin{equation}
v(P) = \min \limits_{x\in \{ 0,1\}^n} \{x^TQx | Ax \le b\} \ge 
\min \limits_{x\in \{ 0,1\} ^n} \{x^TQx + \lambda^T(Ax-b)\}
= v(L_{\lambda})\,.\quad\square\,
\end{equation}

A number of techniques to solve (L) exist in the literature; however, finding
this bound is computationally expensive, so looser bounds (for example, the LP
relaxation) are typically used. Note that the problem yields a natural solution
using the UBQP oracle via the outer Lagrangian linearization method. The book
by Li and Sun~\cite{Li:2006} provides background and several details of this
approach in Procedure 3.2.

Recall that (L) can be rewritten as an LP problem in terms of the real variables
$\lambda$ and $\mu$:
\begin{eqnarray}
	\mathrm{(L_{LP})}& \mbox{Maximize} & \mu \nonumber\\
	& \mbox{subject to} & \mu \leq x^TQx + \lambda^T(Ax-b)
	  \quad (\forall x\in \{ 0,1\}^n)\,, \nonumber\\
	&& \lambda \geq 0\,.
\end{eqnarray}
This formulation is difficult to solve directly, as there are an exponential
number of constraints. In particular, there is one linear constraint (cutting
plane) for every binary point $x \in \{0, 1\}^n$. However, the restriction of
the constraints to a much smaller nonempty subset $T \subseteq \{0, 1\}^n$ of
binary vectors is a tractable LP problem:
\begin{eqnarray}
	\mathrm{(L_{LP-T})}& \mbox{Maximize} & \mu \nonumber\\
	& \mbox{subject to} & \mu \leq x^TQx + \lambda^T(Ax-b)
	  \quad (\forall x\in T)\,, \nonumber\\
	&& \lambda \geq 0\,.
\end{eqnarray}
This LP problem is bounded provided that $T$ contains at least one feasible
solution. In some cases, finding at least one feasible solution can be very
difficult. In these cases, we impose an upper bound on the vector of Lagrange
multipliers, substituting the constraint $\lambda \geq 0$ by $0 \leq \lambda
\leq u$.
\begin{eqnarray}
	\mathrm{(L_{LP-T}^u)}& \mbox{Maximize} & \mu \nonumber\\
	& \mbox{subject to} & \mu \leq x^TQx + \lambda^T(Ax-b)
	  \quad (\forall x\in T)\,, \nonumber\\
	&& 0 \leq \lambda \leq u
\end{eqnarray}

The vector $u$ of upper bounds may depend on an estimation of the solution to
(L). In practical situations, it should also depend on the specifics of the UBQP
oracle (for example, the noise and precision of the oracle). Let $(\mu^\ast,
\lambda^\ast)$ be an optimal solution to $\mathrm{(L_{LP-T})}$. Note that
imposing the box constraint $0 \leq \lambda \leq u$ might result in not
generating an optimal primal-dual pair, but nevertheless $\mathrm{(L_{LP-T}^u)}$
generates a lower bound for $\mathrm{(P)}$. It is clear that
$(\mathrm{L_{\lambda^\ast}})$ is a UBQP problem that can be solved using the
UBQP oracle. By successively adding the solutions returned by the UBQP oracle as
cutting planes and applying the simplex method, we are able to solve
$\mathrm{(L_{LP})}$ (see Algorithm \ref{alg:oa-ld}).

\begin{algorithm}
\caption{Outer Lagrangian linearization}
\label{alg:oa-ld} 
\begin{algorithmic}[1]
	\Input problem $(P)$; initial Lagrange multipliers $\lambda^*$
    \Output optimal primal $x^*$; optimal dual solution $(\lambda^*,\mu^*)$;
    	a strong-duality flag
	\State initialize $\mu^* = \infty$, $T= \emptyset$
	\State solve $\mathrm{L_{\lambda^*}}$, obtain the dual value $d(\lambda^*)$
	and optimal solution $x^* \subseteq \{0, 1\}^n$
	\If{$ (\lambda^*)^T[A(x^*) - b] = 0 $ and $A(x^*) \le b$}
		\State stop and return $x^*$ and $(\lambda^*, d(\lambda^*))$, with
		strong-duality flag set to {\tt true}
	\EndIf
	\If{$ (\mu^*) \leq d(\lambda^*)$}
		\State stop and return $x^*$ and $(\lambda^*, \mu^*)$
	\EndIf
		\State update $T$ by adding $x^*$, $T \leftarrow T \cup \{x^*\}$
		\State solve $L_{LP-T}^u$ and update $\mu^*$ and $\lambda^*$
		\State go to step 2
\end{algorithmic}
\end{algorithm}

\section{A local search heuristic} \label{sec:UB}

In order to prune a branch-and-bound tree effectively, it is important to
quickly obtain a good upper bound. We employ an adaptation of the local search
algorithm presented by Bertsimas et al. \cite{Bertsimas:2013}. The main idea is
as follows: beginning with a feasible solution $x$, the solution is iteratively
improved by considering solutions in the 1-flip neighbourhood of $x$ (defined as
the set of solutions that can be obtained by flipping a single element of $x$)
which are feasible, together with ``interesting'' solutions, until it cannot be
improved further. The algorithm takes a parameter $\rho$ which explicitly
controls the trade-off between complexity and performance by increasing the size
of the neighbourhood. A neighbouring solution $y$ is considered ``interesting''
if it satisfies the following conditions: (i) no constraint is violated by more
than one unit; and (ii) the number of violated constraints in $y$ plus the
number of loose constraints which differ from the loose constraints in the
current best solution is at most $\rho$.

\begin{algorithm}[h]
\caption{Local search heuristic (LSH)}
\label{alg:lsh}
\begin{algorithmic}[1]
	\Input matrices $A$ and $Q$; vector $b$; feasible solution $z_0$; scalar
	  parameter $\rho>0$
    \Output feasible solution $z$ such that $z^TQz\leq z_0^TQz_0$
	\State $z:=z_0$; $S:=\{ z\}$
	\While {$S\neq \emptyset$}
		\State get a new solution $x$ from $S$
		\ForAll {$y$ adjacent to $x$}
			\If{$y$ is feasible and $y^TQy<z^TQz$}
				\State $z\gets y$ and $S\gets y$
				\State go to step 3
			\ElsIf{$y$ is interesting}
				\State $S\gets S\cup \{ y\}$
			\EndIf
		\EndFor
	\EndWhile
\end{algorithmic}
\end{algorithm}

Note that the algorithm moves to a better solution as soon as it finds a
feasible one, and only when no solutions are found does it consider moving to
``interesting'' solutions.

\section{Branching strategies} \label{sec:Branching}
In any branch-and-bound scheme, the performance of the algorithm is largely
dependent on the number of nodes that are visited in the search tree. As such,
it is important to make effective branching decisions, reducing the size of the
search tree. Branching heuristics are usually classified as either static
variable-ordering (SVO) heuristics or dynamic variable-ordering (DVO)
heuristics. All branching heuristics used in this paper are DVO heuristics, as
they are generally considered more effective because they allow information
obtained during a search to be utilized to guide the search.

It is often quite difficult to find an assignment of values to variables that
satisfies all constraints. This has motivated the study of a variety of
approaches that attempt to exploit the interplay between variable-value
assignments and constraints. Examples include the impact-based heuristics
proposed by Refalo~\cite{Refalo:2004}, the conflict-driven variable-ordering
heuristic proposed by Boussemart et al.~\cite{Boussemart:2004}, and the
approximated counting-based heuristics proposed by Kask et al.~\cite{Kask:2004},
Hsu et al.~\cite{Hsu:2007}, Bras et al.~\cite{Bras:2009}, and Pesant et al.
\cite{Pesant:2012}.

\subsection{Counting the solution density} 
\label{sec:maxsd}

One branching heuristic used in this paper is a modified implementation of the
\emph{maxSD} heuristic introduced by Pesant et al.~\cite{Pesant:2012}. We
recall two definitions.

\begin{definition}
	Given a constraint $c(x_1,\dots,x_n)$ and respective finite domains $D_i$
	$(1\leq i\leq n)$, let $\#c(x_1,\ldots ,x_n)$ denote the number of
	$n$-tuples in the corresponding relation.
\end{definition}

\begin{definition}
	Given a constraint $c(x_1,\dots,x_n)$, respective finite domains $D_i$
	$(1\leq i \leq n)$, a variable $x_i$ in the scope of $c$, and a value $d\in
	D_i$, the \emph{solution density} of a pair $(x_i,d)$ in $c$ is given by
	\begin{equation}
		\sigma(x_i,d,c) = \frac{\#c(x_1,\ldots ,x_{i-1},d,x_{i+1},\ldots
		,x_n)}{\#c(x_1,\ldots ,x_n)}\,.
	\end{equation}	
\end{definition}

The solution density measures how frequently a certain assignment of a value in
the domain of a variable belongs to a solution that satisfies constraint $c$.

The heuristic maxSD iterates over all of the variable-value pairs and chooses
the pair that has the highest solution density. If the (approximate)
$\sigma(x_i,d,c)$ are precomputed, the complexity of the modified algorithm is
$O(mq)$, where $m$ is the number of constraints, and $q$ is the sum of the
number of variables that appear in each constraint. Pesant et
al.~\cite{Pesant:2012} detail good approximations of the solution densities for
knapsack constraints, which can be computed efficiently. They also provide an
in-depth experimental analysis that shows that this heuristic is state of the
art among counting-based heuristics.

\subsection{Constraint satisfaction via an LD solution}
\label{sec:cs-bound}

In each node $u$ of the branch-and-bound tree, a lower bound is computed by
solving the LD problem $\mathrm{(L)}$, and the primal-dual pair $(x^u,
\lambda^u)$ is obtained. In the standard way, we define the slack of constraint
$i$ at a point $x$ as $s_i=b_i - a_i^Tx$, where $a_i$ is the $i$-th row of $A$.
Then the set of violated constraints at $x$ is the set $V = \{ i:s_i < 0 \}$. If
$x^u$ is infeasible for the original problem, it must violate one or more
constraints. Additionally, we define the change in slack for constraint $i$
resulting from flipping variable $j$ in $x^u$ as
\begin{equation}
	\delta_{ij} = a_{ij} (2 x^u_{j} -1)\,.
\end{equation}
We present two branching strategies which use this information at $x^u$ to guide
variable and value selection towards feasibility.

The first branching method we propose is to select the variable that maximizes
the reduction in violation of the most violated constraint. That is, we select
$j = \argmax \limits_{j \in 1,\dots n} \delta_{ij}$ and value $1-x^u_j$ (see
Algorithm \ref{alg:scmvb}).
\begin{algorithm}[H]
\caption{Most-violated-constraint satisfaction branching scheme}
\label{alg:scmvb} 
\begin{algorithmic}[1] 
	\Input $x^u$ is the optimal solution to (L) at the current node.
    \ForAll {constraints $i$}
    	\State {compute $s_i = b_i - a^T_ix^u$}
    \EndFor
    \State {$i = \argmin \limits_{i\in \{1,\ldots ,m\}} s_i$}
    \If {$s_i >0$} 
        \State {LD optimal is feasible; abort violation branching} 
    \EndIf
    \ForAll {variables $j$}
    	\State compute $\delta_{ij} = a_{ij} (2 x^u_{j} -1)$ 
    \EndFor
    \State \Return
      index {$j = \argmax \limits_{j \in 1,\dots, n} \delta_{ij}$}
      and value $1 - x^u_j$
\end{algorithmic}
\end{algorithm}

The next branching method we discuss is more general: instead of looking only at
the most violated constraint, we consider all of the violated constraints and
select the variable which, when flipped in the LD solution, gives the maximum
decrease in the left-hand side of all violated constraints (see Algorithm
\ref{alg:allmvb}).

\begin{algorithm}[H]
\caption{All-violated-constraints satisfaction branching scheme}
\label{alg:allmvb} 
\begin{algorithmic}[1] 
	\Input $x^u$ is the optimal solution to (L) at the current node.
    \ForAll {constraints $i$}
    	\State {compute $s_i = b_i - a^T_ix^u$}
    \EndFor
    \State {$V:=\{ i:s_i < 0 \}$ defines the set of violated constraints.}
    \If {$V = \emptyset$} 
        \State {LD optimal is feasible; abort violation branching} 
    \EndIf
    \ForAll {variables $j$}
    	\State compute $\delta_{ij} = a_{ij} (2 x^u_{j} -1)$ 
    \EndFor
    \State \Return 
      index {$j = \argmax \limits_{j \in 1,\dots, n}
      \left(\sum_{i \in V} \delta_{ij}\right)$} and value $1 - x^u_j$
\end{algorithmic}
\end{algorithm}

\subsection{Pseudo-cost branching}

\label{sec:lp-bound}

Introduced in CPLEX 7.5~\cite{Applegate:1995}, the idea of \emph{strong
branching} is to test which fractional variable gives the best bound before
branching. The test is performed by temporarily fixing each fractional variable
to 0 or 1 and solving the LP relaxation by the dual simplex method. Since the
cost of solving several LP subproblems is high, only a fixed number of
iterations of the dual simplex algorithm are performed. The variable fixation
that provides the strongest bound is chosen as the branching decision.

If the number of fractional variables is large, this process is very time
consuming. There are several possible methods to overcome this difficulty. One
way is to select a subset of variables, for example, choosing variables with
values close to 0.5. Another approach, the \emph{$k$-look-ahead} branching
strategy, requires an integer parameter $k$ and a score assignment on the
fractional variable fixations. The variable fixations are then sorted according
to their scores, and the above test is performed. If no better bound is found
for $k$ successive variable fixations, the test process is stopped.

The sorting of the variable fixations is only applied in later stages of the
branch-and-bound process. We use \emph{pseudo-costs}, which are introduced in
\cite{pseudocost} and explained in Section \ref{sec:pseudo-cost-def}, as the
score of a variable fixation.

Note that the lower-bound computation by the UBQP oracle can be performed in
parallel with either strong branching or the $k$-look-ahead strategy using a
digital processor, affording the computational time required to utilize this
type of branching without significantly increasing the total running time of the
algorithm.

\subsubsection{Pseudo-costs}\label{sec:pseudo-cost-def}
Pseudo-costs keep track of the success of the variable fixations that have
already been used in the branch-and-bound process. Different variations of
pseudo-costs have been proposed in the literature. We employ the variation
discussed in \cite{Achterberg,Martin1998}. Let $f_i^+$ and $f_i^-$ be the amount
that the variable $x_i$ is respectively rounded up and down:
\begin{equation}
f_i^+ = \lceil x_i \rceil -x_i\,, \quad
\text{and} \quad f_i^- = x_i - \lfloor x_i \rfloor\,.
\end{equation}
We let $(\mathrm S)$ denote the subproblem associated to a node in the branch-
and-bound tree that has already been visited and branched on. Let $i$ denote the
index of the variable $x_i$ over which the algorithm branched from $(\mathrm
S)$. We use the notation $\mathrm{S_i^+}$ and $\mathrm{S_i^-}$ for the child
subproblems of $(\mathrm S)$ by branching on the variable $x_i$ to $1$ and $0$,
respectively. We define
\begin{equation}
\Delta_{i}^+ = v(\mathrm{(S_i^+)_{LP}}) - v(\mathrm{S_{LP}})\,, 
\quad\text{ and }\quad
\Delta_{i}^- = v(\mathrm{(S_i^-)_{LP}}) - v(\mathrm{S_{LP}})\,,
\end{equation}
as the change in the optimal values of the linear relaxations of problems
$\mathrm{(S_i^+)}$ and $\mathrm{(S_i^-)}$ from that of $\mathrm{(S)}$. Let
$\mathcal{G}_i^+$ and $\mathcal{G}_i^-$ be the above difference per unit of
change in variable $x_i$ at subproblem $\mathrm{(S)}$:
\begin{equation} 
\mathcal{G}_i^+ = \frac{\Delta_{i}^+}{f_i^+}\,,\quad
\text{ and } \quad\mathcal{G}_i^- = \frac{\Delta_{i}^-}{f_i^-}\,.
\end{equation}
We let $\delta_i^+$ be the sum of $\mathcal{G}_i^+$ over all subproblems
$\mathrm S$ for which fixation of $x_i$ to $1$ is selected and the LP-relaxation
of the subproblem $\mathrm{(S)_i^+}$ was feasible. $\mathcal N_i^+$ is the
number of all such subproblems. $\delta_i^-$ and $\mathcal N_i^-$ are defined
analogously for the fixation of $x_i$ to $0$.

Then the pseudo-cost of upward and downward branching of variable $x_i$ are
defined as
\begin{equation}
\Upsilon_i^+ = \frac{\delta_i^+}{\mathcal{N}_i^+}\,, \quad \text{ and } \quad
\Upsilon_i^- = \frac{\delta_i^-}{\mathcal{N}_i^-}\,.
\end{equation}
Finally, the score assigned to the variable fixations of $x_i$ to 0 and 1 is the
following convex combination
\begin{equation}
s_{x_i} = (1-\mu). \mathrm{min}(f_i^+ \Upsilon_i^+,f_i^- \Upsilon_i^-) + \mu.
\mathrm{max}(f_i^+ \Upsilon_i^+,f_i^- \Upsilon_i^-)\,.
\end{equation}
The score factor $\mu$ is a number between 0 and 1. In our experiments, this
factor is set to 0.3.

\subsection{Frequency-based branching} 
\label{sec:freq-bound}

Motivated by the notion of \emph{persistencies} as described by Boros and
Hammer~\cite{Boros:2002}, and observing that the outer Lagrangian linearization
method yields a number of high-quality solutions, one can perform $k$-look-ahead
branching, selecting variable-value pairs based on their \emph{frequency}. Here,
given a set $S \subseteq \{0, 1\}^n$ of binary vectors, an index $i\in \{1,
\ldots, n\}$, and a binary value $s \in \{0, 1\}$, the frequency of the pair
$(x_i, s)$ in $S$ is defined as the number of elements in $S$ with their $i$-th
entry equal to $s$. When using a UBQP oracle that performs quantum annealing,
the oracle returns a spectrum of solutions and all solutions can be used in the
frequency calculation. This branching strategy has not, to our knowledge,
previously appeared in the literature.

\section{Branch-and-bound framework}
\label{sec:qbnb}

We now present our branch-and-bound algorithm in its entirety, before reporting
the computational results of its performance. The computation of the Lagrangian
dual bound is skipped at every node unless a finite upper bound exists, that is,
a feasible solution is known. If a feasible solution is not yet observed, the
maxSD heuristic of Section \ref{sec:maxsd} is used for branching. Once a
feasible solution is observed, the Lagrangian dual bounds are computed and the
branching strategy switches from maxSD to one of the bounding methods explained
in Sections \ref{sec:cs-bound}, \ref{sec:lp-bound}, and \ref{sec:freq-bound}.

After the first feasible solution is found, the heuristic of Algorithm
\ref{alg:lsh} is executed on another processor core and improves the best upper
bound found thus far in parallel to the branch-and-bound algorithm.

\begin{algorithm}[h]
\caption{Branch and bound}
\begin{algorithmic}[1]
	\Input matrices $A$ and $Q$; vector $b$
    	\Output optimal solution value $z^*$; or report there is no optimal solution
	\State $z = \infty$ , $k = 0$, $\rho =1$ , $L= \{(P)\}$
	\While {$L \neq \emptyset$} 
	\State choose a problem $p$ from $L$
	\If {$p$ has no variables} 
		\State set $c$ to the constant objective of $p$
		\If {the fixed variables create a feasible solution and $c < z$}
		\State $z \leftarrow c$ and update optimal solution
		\EndIf
	\ElsIf {$z = \infty$} 
		\State create subproblems $p_1$ and $p_2$ according to maxSD
		\State update $L \leftarrow L \cup \{p_1, p_2\}$
	\Else
		\State perform outer Lagrangian linearization on $p$ 
		\State obtain $x^*$, $(\lambda^*, \mu^*)$ and strong duality flag
		\If {$\mu^* \leq z$}
			\If {$x^*$ and fixed variables create a feasible solution and $\mu^* < z$} 
				\State $z \leftarrow \mu^*$ and update optimal solution
			\EndIf
			\If {strong duality is {\tt false}} 
				\State create subproblems $p_1$ and $p_2$ according to a branching strategy
				\State update $L \leftarrow L \cup \{p_1, p_2\}$
			\EndIf
		\EndIf
	\EndIf
	\EndWhile
     \end{algorithmic}
\end{algorithm}

\section{Computational experiments}\label{sec:Experiments}

\subsection{Generation of test instances}
For this paper, we used randomly generated test instances with inequality
constraints. For a specific test instance, we let $n$ denote the number of
variables and $m$ be the number of inequality constraints. To generate the cost
matrix, we constructed an $n \times n$ random symmetric matrix $Q$ of a given
density $d$. An $m \times n$ constraint matrix was generated in a similar
manner, ensuring that the CBQP problem had at least one feasible solution.
Densities of 0.3 and 0.5 were used for the objective functions and the
constraint matrices, respectively. The values of $n$ ranged between $36$ and
$50$, and $10$ instances were generated for each size. The values of $m$ were
chosen as $\frac{n}2$ for even $n$. When $n$ was odd, $5$ of the instances had $m
= \frac{n-1}2$ and the other $5$ instances had $m= \frac{n+1}2$ as numbers of
constraints.

\subsection{Computational results}

We now present the details of our computational experiments. The algorithms
were programmed in C++ and compiled using GNU~GCC on a machine with a 2.5~GHz
Intel Core~i5-3210M processor and 16~GB of RAM. Linear programming and UBQP
problems were solved by the Gurobi Optimizer~5.6, using the Gurobi Optimizer in
place of a UBQP oracle and replacing the computational time with 0~milliseconds
per solver call.  The algorithm was coded to utilize 4 cores and to allow us to
accurately report times. All other threads were paused during the solving of the
UBQP problems.

\newlength{\defaulttabcolsep}
\setlength{\defaulttabcolsep}{\tabcolsep}
\setlength{\tabcolsep}{1pt}
\renewcommand{\arraystretch}{.65}

\begin{table}[b!]
\caption{Branch-and-bound nodes needed to solve each instance.}
\label{tab:t1}
\makeatletter\@setfontsize\minisclue{8}{10}\makeatother
	\begin{tabular}{rl}
		\begin{tabular}[t]{ccccccccc}
		\toprule
			{n} & {m} & {mviol} & {aviol} & {pc4} & {pc8} & {freq4} & {freq8} & {gurobi}\\
			\midrule
			$36$ & $18$ & $35$ & $35$ & $41$ & $37$ & $53$ & \textbf{15} & $4336$\\
			$36$ & $18$ & $59$ & $67$ & $57$ & \textbf{39} & $315$ & 121 & $2699$\\
			$36$ & $18$ & $27$ & $27$ & $27$ & \textbf{19} & $45$ & 45 & $2411$\\
			$36$ & $18$ & \textbf{7} & \textbf{7} & $67$ & 57 & $103$ & 13 & $6381$\\
			$36$ & $18$ & \textbf{3} & \textbf{3} & $39$ & 31 & $5$ & 5 & $12003$\\
			$36$ & $18$ & \textbf{3} & \textbf{3} & $37$ & 25 & $11$ & 5 & $11459$\\
			$36$ & $18$ & \textbf{11} & 13 & $49$ & 23 & $51$ & 23 & $18004$\\
			$36$ & $18$ & \textbf{5} & 7 & \textbf{5} & \textbf{5} & $27$ & 25 & $564$\\
			$36$ & $18$ & 133 & \textbf{109} & 175 & 169 & $207$ & 139 & $28216$\\
			$36$ & $18$ & 21 & 21 & 61 & 31 & $23$ & \textbf{13} & $6311$\\
			$37$ & $18$ & \textbf{7} & \textbf{7} & 13 & 13 & $31$ & 41 & $2512$\\
			$37$ & $18$ & 27 & \textbf{23} & 69 & 77 & $31$ & 29 & $20452$\\
			$37$ & $18$ & \textbf{13} & 17 & 87 & 63 & $161$ & 93 & $5553$\\
			$37$ & $18$ & \textbf{3} & \textbf{3} & 7 & \textbf{3} & \textbf{3} & \textbf{3} & $20642$\\
			$37$ & $18$ & 13 & 7 & 81 & 71 & \textbf{5} & \textbf{5} & $13490$\\
			$37$ & $19$ & 49 & 49 & 51 & 55 & 53 & \textbf{37} & $18079$\\
			$37$ & $19$ & 43 & 43 & \textbf{19} & \textbf{19} & 67 & 43 & $5644$\\
			$37$ & $19$ & 47 & \textbf{33} & 65 & 49 & 71 & 37 & $17504$\\
			$37$ & $19$ & 9 & \textbf{5} & 55 & 55 & 73 & 13 & $26857$\\
			$37$ & $19$ & 13 & 11 & 17 & 17 & \textbf{5} & \textbf{5} & $16055$\\
			$38$ & $19$ & 43 & \textbf{41} & 113 & 113 & 201 & 177 & $3754$\\
			$38$ & $19$ & \textbf{5} & \textbf{5} & 75 & 43 & 29 & 37 & $2572$\\
			$38$ & $19$ & \textbf{5} & \textbf{5} & 51 & 47 & 69 & 45 & $13306$\\
			$38$ & $19$ & 11 & \textbf{7} & 165 & 103 & 105 & 75 & $16847$\\
			$38$ & $19$ & 13 & 13 & 51 & 29 & 11 & \textbf{9} & $22237$\\
			$38$ & $19$ & \textbf{5} & \textbf{5} & 73 & 57 & 51 & 11 & $13366$\\
			$38$ & $19$ & \textbf{3} & \textbf{3} & 11 & 11 & 9 & 7 & $17339$\\
			$38$ & $19$ & 71 & \textbf{47} & 69 & 75 & 69 & 75 & $67932$\\
			$38$ & $19$ & \textbf{5} & \textbf{5} & 81 & 55 & 95 & 45 & $13313$\\
			$38$ & $19$ & 23 & 23 & 113 & 63 & \textbf{21} & 41 & $16600$\\
			$39$ & $19$ & 19 & \textbf{17} & 79 & 57 & 31 & 21 & $6019$\\
			$39$ & $19$ & 123 & \textbf{97} & 231 & 169 & 143 & 135 & $64884$\\
			$39$ & $19$ & \textbf{37} & \textbf{37} & 131 & 77 & 253 & 51 & $14214$\\
			$39$ & $19$ & \textbf{25} & \textbf{25} & 147 & 85 & 127 & 131 & $11538$\\
			$39$ & $19$ & \textbf{9} & \textbf{9} & 41 & 21 & 55 & 29 & $7175$\\
			$39$ & $20$ & 9 & 17 & 65 & 65 & 23 & \textbf{7} & $75381$\\
			$39$ & $20$ & 45 & 45 & \textbf{21} & 39 & 75 & 59 & $5739$\\
			$39$ & $20$ & 35 & 27 & \textbf{15} & \textbf{15} & 93 & 63 & $1488$\\
			$39$ & $20$ & \textbf{5} & \textbf{5} & 11 & 39 & 9 & 9 & $3518$\\
			$39$ & $20$ & 7 & 7 & 131 & 31 & \textbf{5} & 11 & $88351$\\
			$40$ & $20$ & \textbf{3} & \textbf{3} & 95 & 63 & 79 & 11 & $39004$\\
			$40$ & $20$ & \textbf{3} & \textbf{3} & 23 & 27 & 11 & 11 & $9917$\\
			$40$ & $20$ & \textbf{17} & \textbf{17} & 215 & 103 & 173 & 145 & $13934$\\
			$40$ & $20$ & \textbf{5} & \textbf{5} & 9 & 55 & 29 & 29 & $7482$\\
			$40$ & $20$ & \textbf{23} & 43 & 131 & 59 & 221 & 81 & $4481$\\
			$40$ & $20$ & \textbf{9} & \textbf{9} & 125 & 61 & 53 & 23 & $127674$\\
			$40$ & $20$ & 37 & 21 & 111 & 87 & 119 & \textbf{15} & $19213$\\
			$40$ & $20$ & 19 & 19 & 39 & 31 & 25 & \textbf{17} & $56666$\\
			$40$ & $20$ & 47 & \textbf{43} & 189 & 115 & 93 & 47 & $5617$\\
			$40$ & $20$ & \textbf{25} & \textbf{25} & 179 & 89 & 49 & 115 & $32531$\\
			$41$ & $20$ & \textbf{9} & \textbf{9} & 79 & 47 & 25 & \textbf{9} & $97004$\\
			$41$ & $20$ & \textbf{3} & \textbf{3} & 7 & 19 & 5 & 5 & $48826$\\
			$41$ & $20$ & \textbf{7} & \textbf{7} & 43 & 87 & 17 & 17 & $48518$\\
			$41$ & $20$ & 45 & \textbf{39} & 289 & 87 & 281 & 91 & $28629$\\
			$41$ & $20$ & \textbf{23} & 27 & 97 & 43 & 125 & 59 & $65260$\\
			$41$ & $21$ & \textbf{5} & \textbf{5} & 75 & 29 & 35 & \textbf{5} & $37813$\\
			$41$ & $21$ & \textbf{11} & \textbf{11} & 41 & 31 & 99 & 131 & $8861$\\
			$41$ & $21$ & 15 & 15 & 15 & 33 & 13 & \textbf{11} & $1970$\\
			$41$ & $21$ & \textbf{7} & 11 & 61 & 81 & 51 & 43 & $40794$\\
			$41$ & $21$ & \textbf{45} & 55 & 297 & 213 & 531 & 81 & $20526$\\
			$42$ & $21$ & \textbf{21} & 25 & 67 & 125 & 219 & 45 & $79671$\\
			$42$ & $21$ & 11 & 9 & 51 & 59 & \textbf{5} & 53 & $15149$\\
			$42$ & $21$ & \textbf{103} & \textbf{103} & 123 & 107 & 229 & 177 & $46238$\\
			$42$ & $21$ & \textbf{5} & \textbf{5} & 159 & 155 & \textbf{5} & \textbf{5} & $55574$\\
			$42$ & $21$ & 33 & \textbf{31} & 231 & 97 & 133 & 103 & $8623$\\
			$42$ & $21$ & 31 & \textbf{27} & 341 & 223 & 125 & 69 & $80136$\\
			$42$ & $21$ & 17 & 49 & 45 & 37 & 37 & \textbf{11} & $96842$\\
			$42$ & $21$ & \textbf{31} & \textbf{31} & 61 & 93 & 51 & 47 & $59757$\\
			$42$ & $21$ & \textbf{45} & \textbf{45} & 123 & 99 & 175 & 75 & $3605$\\
			$42$ & $21$ & 35 & 35 & 45 & \textbf{29} & 53 & 87 & $114468$\\
			$43$ & $21$ & \textbf{5} & \textbf{5} & 103 & 55 & 111 & \textbf{5} & $52459$\\
			$43$ & $21$ & \textbf{41} & \textbf{41} & 253 & 149 & 67 & 133 & $512973$\\
			$43$ & $21$ & 37 & 37 & \textbf{5} & 31 & 11 & 7 & $48384$\\
			$43$ & $21$ & \textbf{3} & \textbf{3} & 53 & 55 & 105 & 67 & $15789$\\
			$43$ & $21$ & 15 & 25 & 31 & \textbf{5} & 31 & 13 & $132226$\\
		\bottomrule
		\end{tabular}
		&\quad
		\begin{tabular}[t]{ccccccccc}
		\toprule
			{n} & {m} & {mviol} & {aviol} & {pc4} & {pc8} & {freq4} & {freq8} & {gurobi}\\
			\midrule
			$43$ & $22$ & 77 & \textbf{61} & 407 & 351 & 387 & 257 & $67229$\\
			$43$ & $22$ & \textbf{9} & \textbf{9} & 175 & 179 & 59 & 63 & $31891$\\
			$43$ & $22$ & \textbf{7} & \textbf{7} & 85 & 97 & 109 & 111 & $10037$\\
			$43$ & $22$ & 19 & \textbf{17} & 287 & 185 & 29 & 79 & $110101$\\
			$43$ & $22$ & \textbf{5} & \textbf{5} & 25 & 35 & 7 & 7 & $68284$\\
			$44$ & $22$ & \textbf{7} & \textbf{7} & 89 & 77 & 43 & 23 & $123583$\\
			$44$ & $22$ & \textbf{13} & \textbf{13} & 269 & 127 & 87 & 163 & $190688$\\
			$44$ & $22$ & \textbf{5} & \textbf{5} & 47 & 67 & 77 & 51 & $6865$\\
			$44$ & $22$ & \textbf{85} & 117 & 223 & 209 & 427 & 375 & $169577$\\
			$44$ & $22$ & 143 & 117 & 173 & \textbf{109} & 245 & 191 & $49607$\\
			$44$ & $22$ & 89 & \textbf{73} & 531 & 263 & 317 & 79 & $563781$\\
			$44$ & $22$ & 65 & 73 & 195 & 213 & 137 & \textbf{61} & $805637$\\
			$44$ & $22$ & 363 & 391 & 451 & \textbf{353} & 681 & 469 & $424882$\\
			$44$ & $22$ & \textbf{21} & 31 & 43 & 23 & 27 & 137 & $111867$\\
			$44$ & $22$ & 7 & \textbf{5} & 95 & 97 & 107 & 13 & $11674$\\
			$45$ & $22$ & \textbf{5} & \textbf{5} & 113 & 25 & 73 & 11 & $162930$\\
			$45$ & $22$ & 491 & \textbf{169} & 757 & 545 & 577 & 305 & $476795$\\
			$45$ & $22$ & \textbf{5} & \textbf{5} & 149 & 109 & 9 & 9 & $85671$\\
			$45$ & $22$ & \textbf{5} & 7 & 135 & 99 & 45 & 17 & $110672$\\
			$45$ & $22$ & \textbf{5} & \textbf{5} & 145 & 33 & 13 & \textbf{5} & $90347$\\
			$45$ & $23$ & \textbf{5} & \textbf{5} & 97 & 53 & 41 & 23 & $576684$\\
			$45$ & $23$ & \textbf{41} & \textbf{41} & 343 & 141 & 191 & 215 & $77693$\\
			$45$ & $23$ & \textbf{15} & 23 & 269 & 145 & 143 & 37 & $163173$\\
			$45$ & $23$ & \textbf{87} & 125 & 289 & 231 & 161 & 133 & $10268$\\
			$45$ & $23$ & \textbf{91} & 95 & 215 & 259 & 599 & 173 & $140421$\\
			$46$ & $23$ & 11 & 11 & 27 & 201 & \textbf{5} & \textbf{5} & $120440$\\
			$46$ & $23$ & \textbf{19} & \textbf{19} & 243 & 237 & 131 & 117 & $53680$\\
			$46$ & $23$ & 39 & 39 & 41 & 27 & 39 & \textbf{25} & $243014$\\
			$46$ & $23$ & 113 & \textbf{69} & 441 & 297 & 713 & 311 & $90499$\\
			$46$ & $23$ & 29 & \textbf{21} & 161 & 37 & 37 & 43 & $63781$\\
			$46$ & $23$ & 269 & \textbf{105} & 441 & 331 & 869 & 329 & $134925$\\
			$46$ & $23$ & \textbf{13} & \textbf{13} & 201 & 165 & 145 & 77 & $39956$\\
			$46$ & $23$ & \textbf{13} & \textbf{13} & 163 & 145 & 83 & 49 & $249360$\\
			$46$ & $23$ & \textbf{11} & \textbf{11} & 25 & 75 & 15 & \textbf{11} & $1082407$\\
			$46$ & $23$ & \textbf{7} & \textbf{7} & 249 & 119 & 193 & 69 & $352602$\\
			$47$ & $23$ & \textbf{9} & \textbf{9} & 287 & 181 & 93 & 21 & $2535458$\\
			$47$ & $23$ & \textbf{9} & 13 & 139 & 97 & 115 & 13 & $763493$\\
			$47$ & $23$ & \textbf{25} & 59 & 725 & 419 & 1095 & 455 & $323064$\\
			$47$ & $23$ & 7 & 7 & \textbf{3} & \textbf{3} & 5 & 5 & $526230$\\
			$47$ & $23$ & 41 & 57 & 37 & \textbf{5} & 223 & 167 & $191022$\\
			$47$ & $24$ & 31 & 31 & 167 & 97 & 77 & \textbf{27} & $34062$\\
			$47$ & $24$ & \textbf{55} & \textbf{55} & 75 & 95 & 129 & 137 & $181039$\\
			$47$ & $24$ & \textbf{19} & \textbf{19} & 465 & 365 & 305 & 173 & $459243$\\
			$47$ & $24$ & \textbf{53} & 61 & 683 & 373 & 515 & 377 & $62009$\\
			$47$ & $24$ & \textbf{15} & \textbf{15} & 57 & 173 & 565 & 133 & $187353$\\
			$48$ & $24$ & \textbf{5} & \textbf{5} & 73 & 109 & 29 & 15 & $6302$\\
			$48$ & $24$ & \textbf{33} & \textbf{33} & 125 & 99 & 263 & 81 & $200975$\\
			$48$ & $24$ & \textbf{61} & \textbf{61} & 175 & 155 & 293 & 235 & $304816$\\
			$48$ & $24$ & \textbf{5} & \textbf{5} & 81 & 79 & 137 & 17 & $214248$\\
			$48$ & $24$ & 59 & \textbf{49} & 211 & 271 & 561 & 171 & $71703$\\
			$48$ & $24$ & \textbf{31} & 33 & 91 & 215 & 209 & 97 & $47419$\\
			$48$ & $24$ & \textbf{15} & \textbf{15} & 121 & 131 & 31 & 19 & $2977412$\\
			$48$ & $24$ & 23 & 23 & 175 & 117 & 23 & \textbf{21} & $3685777$\\
			$48$ & $24$ & \textbf{9} & 13 & 421 & 225 & 469 & 31 & $347306$\\
			$48$ & $24$ & \textbf{11} & \textbf{11} & 83 & 69 & 31 & 71 & $655210$\\
			$49$ & $24$ & \textbf{83} & 85 & 841 & 469 & 297 & 337 & $244676$\\
			$49$ & $24$ & 27 & 19 & 95 & 75 & 63 & \textbf{7} & $7137848$\\
			$49$ & $24$ & 27 & \textbf{19} & 279 & 139 & 109 & 57 & $7261850$\\
			$49$ & $24$ & \textbf{19} & \textbf{19} & 319 & 241 & 39 & 25 & $898477$\\
			$49$ & $24$ & \textbf{15} & \textbf{15} & 131 & 79 & 29 & 41 & $524699$\\
			$49$ & $25$ & \textbf{11} & 21 & 421 & 197 & 301 & 31 & $5346494$\\
			$49$ & $25$ & 57 & \textbf{55} & 215 & 203 & 337 & 369 & $439194$\\
			$49$ & $25$ & 65 & \textbf{31} & 253 & 203 & 73 & 47 & $1092470$\\
			$49$ & $25$ & 59 & \textbf{35} & 253 & 195 & 107 & 65 & $1014756$\\
			$49$ & $25$ & 9 & 9 & 353 & \textbf{3} & 19 & 23 & $2586700$\\
			$50$ & $25$ & \textbf{33} & 37 & 275 & 209 & 131 & 91 & $3286299$\\
			$50$ & $25$ & 7 & 7 & \textbf{5} & 279 & 61 & 41 & $5763159$\\
			$50$ & $25$ & \textbf{75} & \textbf{75} & 303 & 149 & 139 & 95 & $2763317$\\
			$50$ & $25$ & \textbf{27} & 39 & 331 & 257 & 207 & 73 & $4654692$\\
			$50$ & $25$ & \textbf{31} & \textbf{31} & 297 & 301 & 387 & 93 & $650421$\\
			$50$ & $25$ & \textbf{11} & \textbf{11} & 201 & 111 & 135 & 49 & $2249101$\\
			$50$ & $25$ & \textbf{7} & \textbf{7} & 51 & 55 & 101 & 45 & $230022$\\
			$50$ & $25$ & \textbf{23} & \textbf{23} & 83 & 345 & 399 & 319 & $322977$\\
			$50$ & $25$ & 19 & 33 & 47 & 35 & \textbf{15} & 19 & $415582$\\
			$50$ & $25$ & 37 & \textbf{27} & 331 & 261 & 143 & 75 & $1729174$\\
		\bottomrule
		\end{tabular}
	\end{tabular}
\end{table}

\renewcommand{\arraystretch}{.64}
\begin{table}[b!]
\caption{Time required to solve each instance.}
\label{tab:t2}
\makeatletter\@setfontsize\minisclue{9}{8}\makeatother
	\begin{tabular}{rl}
		\begin{tabular}[t]{ccccccccccc}
		\toprule
			{n} & {m} & {mviol} & {aviol} & {pc4} & {pc8} & {freq4} & {freq8} & {gurobi} & {nq} & {qal}\\
			\midrule
			$36$ & $18$ & $1.67$ & $1.73$ & $0.56$ & $0.96$ & $0.69$ & $0.7$ & $0.23$ & $75$ & $-4$\\
			$36$ & $18$ & $2.29$ & $2.91$ & $1.35$ & $1.57$ & $2.86$ & $1.38$ & $0.2$ & $135$ & $-9$\\
			$36$ & $18$ & $1.57$ & $1.64$ & $0.73$ & $1.07$ & $0.67$ & $0.92$ & $0.18$ & $65$ & $-8$\\
			$36$ & $18$ & $0.7$ & $0.7$ & $1.62$ & $1.47$ & $1.07$ & $0.59$ & $0.38$ & $27$ & $-8$\\
			$36$ & $18$ & $0.34$ & $0.34$ & $0.58$ & $0.79$ & \textbf{0.16} & $0.19$ & $0.67$ & $19$ & $27$\\
			$36$ & $18$ & $0.41$ & $0.33$ & $0.41$ & $0.52$ & 0.31 & \textbf{0.25} & $0.6$ & $13$ & $27$\\
			$36$ & $18$ & $0.65$ & \textbf{0.46} & $0.68$ & $0.77$ & 0.76 & 1.02 & $1.06$ & $36$ & $17$\\
			$36$ & $18$ & $0.52$ & 0.58 & $0.59$ & $0.54$ & 0.38 & 0.37 & $0.05$ & $32$ & $-10$\\
			$36$ & $18$ & $4.11$ & 3.89 & $2.34$ & $2.89$ & 3.37 & 4.26 & $1.6$ & $660$ & $-1$\\
			$36$ & $18$ & $1.9$ & 1.77 & $1$ & $0.86$ & 0.39 & 0.46 & $0.38$ & $73$ & $-0$\\
			$37$ & $18$ & $0.72$ & 0.77 & $1.31$ & $1.23$ & 0.72 & 0.86 & $0.18$ & $45$ & $-12$\\
			$37$ & $18$ & $1.76$ & 1.25 & $1.12$ & $1.5$ & \textbf{0.67} & 1.09 & $1.27$ & $113$ & $5$\\
			$37$ & $18$ & $0.8$ & 1.04 & $1.41$ & $1.68$ & 1.63 & 1.21 & $0.39$ & $73$ & $-6$\\
			$37$ & $18$ & $0.36$ & 0.36 & $0.1$ & \textbf{0.09} & 0.1 & 0.11 & $1.05$ & $11$ & $87$\\
			$37$ & $18$ & $1.38$ & 0.8 & $1.74$ & 1.71 & \textbf{0.2} & 0.32 & $0.73$ & $18$ & $29$\\
			$37$ & $19$ & $1.92$ & 1.94 & \textbf{0.54} & 1.21 & 0.71 & 0.93 & $1$ & $119$ & $4$\\
			$37$ & $19$ & $3.29$ & 3.39 & 0.6 & 0.66 & 0.69 & 0.69 & $0.35$ & $30$ & $-8$\\
			$37$ & $19$ & $1.69$ & 1.35 & \textbf{0.93} & 1.14 & 1.33 & 1.42 & $1.01$ & $148$ & $1$\\
			$37$ & $19$ & $0.44$ & \textbf{0.4} & 0.67 & 1.05 & 0.9 & 0.53 & $1.6$ & $31$ & $39$\\
			$37$ & $19$ & $0.6$ & 0.37 & 0.27 & 0.42 & \textbf{0.17} & 0.19 & $0.93$ & $17$ & $45$\\
			$38$ & $19$ & $2.38$ & 2.26 & 1.86 & 1.76 & 2.16 & 2.21 & $0.34$ & $373$ & $-4$\\
			$38$ & $19$ & $0.64$ & 0.65 & 1.5 & 1.28 & 0.98 & 0.92 & $0.19$ & $32$ & $-14$\\
			$38$ & $19$ & $0.77$ & 0.78 & 0.84 & 1.04 & 0.9 & \textbf{0.66} & $0.77$ & $73$ & $2$\\
			$38$ & $19$ & $1.15$ & \textbf{1.08} & 3.37 & 2.46 & 2.98 & 1.8 & $1.15$ & $77$ & $1$\\
			$38$ & $19$ & $0.53$ & 0.53 & 1.3 & 0.8 & 0.52 & \textbf{0.48} & $1.51$ & $37$ & $28$\\
			$38$ & $19$ & $0.7$ & 0.69 & 1.47 & 1.64 & 0.57 & \textbf{0.42} & $0.83$ & $54$ & $8$\\
			$38$ & $19$ & $0.36$ & 0.37 & 0.28 & \textbf{0.24} & 0.39 & 0.48 & $0.94$ & $33$ & $21$\\
			$38$ & $19$ & $2.42$ & 1.81 & \textbf{0.96} & 1.19 & 1.39 & 1.75 & $3.67$ & $182$ & $15$\\
			$38$ & $19$ & \textbf{0.79} & \textbf{0.79} & 2.03 & 2.03 & 1.2 & 1.17 & $0.94$ & $22$ & $7$\\
			$38$ & $19$ & \textbf{0.89} & 1.29 & 2.16 & 1.68 & 1.13 & 2.31 & $1.38$ & $205$ & $2$\\
			$39$ & $19$ & 0.79 & 0.75 & 1.32 & 1.44 & 1.24 & 1.27 & $0.45$ & $130$ & $-2$\\
			$39$ & $19$ & \textbf{2.06} & 3.56 & 3.72 & 3.42 & 2.45 & 3.02 & $3.99$ & $482$ & $4$\\
			$39$ & $19$ & 2.42 & 2.61 & 2.5 & 2.45 & 5.3 & 2.83 & $0.95$ & $87$ & $-17$\\
			$39$ & $19$ & 3.23 & 3.25 & 4.2 & 3.16 & 1.53 & 1.66 & $0.85$ & $155$ & $-4$\\
			$39$ & $19$ & 1.8 & 1.8 & 3.03 & 1.66 & 0.82 & 0.72 & $0.58$ & $64$ & $-2$\\
			$39$ & $20$ & 0.94 & 0.94 & 0.76 & 1.19 & 0.79 & \textbf{0.46} & $4.61$ & $30$ & $138$\\
			$39$ & $20$ & 2.11 & 2.1 & \textbf{0.31} & 0.74 & 0.66 & 0.65 & $0.38$ & $34$ & $2$\\
			$39$ & $20$ & 2.23 & 2.14 & 1.02 & 1.0 & 0.78 & 0.76 & $0.12$ & $58$ & $-11$\\
			$39$ & $20$ & 1.02 & 1.01 & \textbf{0.21} & 1.29 & 0.24 & 0.26 & $0.26$ & $20$ & $3$\\
			$39$ & $20$ & 0.77 & 0.76 & 2.27 & 0.27 & \textbf{0.22} & 0.57 & $5.08$ & $29$ & $168$\\
			$40$ & $20$ & 0.76 & 0.76 & 2.23 & 1.81 & 2.05 & \textbf{0.66} & $2.62$ & $42$ & $47$\\
			$40$ & $20$ & 0.93 & 0.9 & 0.95 & 1.09 & 0.68 & 0.92 & $0.49$ & $17$ & $-11$\\
			$40$ & $20$ & 1.22 & 1.22 & 4.54 & 3.02 & 1.68 & 3.76 & $0.9$ & $63$ & $-5$\\
			$40$ & $20$ & 0.38 & \textbf{0.37} & 0.81 & 2.43 & 1.02 & 1.16 & $0.54$ & $11$ & $15$\\
			$40$ & $20$ & 1.94 & 1.24 & 3.13 & 2.32 & 2.08 & 1.84 & $0.46$ & $74$ & $-11$\\
			$40$ & $20$ & \textbf{0.82} & \textbf{0.82} & 1.6 & 1.03 & 1.24 & 1.27 & $7.37$ & $31$ & $211$\\
			$40$ & $20$ & 2.87 & 1.91 & 3.1 & 6.41 & 3.45 & \textbf{1.27} & $1.35$ & $49$ & $2$\\
			$40$ & $20$ & 1.48 & 1.47 & \textbf{0.57} & 0.79 & 0.95 & 0.98 & $3.47$ & $230$ & $13$\\
			$40$ & $20$ & 2.61 & 2.81 & 3.24 & 3.34 & 2.95 & 2.99 & $0.4$ & $196$ & $-11$\\
			$40$ & $20$ & 1.74 & 1.69 & 4.17 & 2.84 & \textbf{1.24} & 2.19 & $2.09$ & $121$ & $7$\\
			$41$ & $20$ & \textbf{0.53} & 0.54 & 1.18 & 1.03 & 0.59 & 0.55 & $5.69$ & $32$ & $161$\\
			$41$ & $20$ & 0.44 & 0.45 & \textbf{0.12} & 0.47 & 0.17 & 0.25 & $2.82$ & $15$ & $180$\\
			$41$ & $20$ & 0.5 & \textbf{0.48} & 0.57 & 2.24 & 0.69 & 0.93 & $3.04$ & $16$ & $160$\\
			$41$ & $20$ & 2.4 & 2.48 & 10.91 & 4.46 & 4.0 & 2.65 & $1.78$ & $166$ & $-4$\\
			$41$ & $20$ & \textbf{1.04} & 2.04 & 2.2 & 1.1 & 1.9 & 1.46 & $4.03$ & $103$ & $29$\\
			$41$ & $21$ & 0.97 & 0.97 & 1.69 & 1.24 & 0.72 & \textbf{0.31} & $2.66$ & $17$ & $138$\\
			$41$ & $21$ & 1.64 & 1.65 & \textbf{0.55} & 1.47 & 1.84 & 2.1 & $0.77$ & $102$ & $2$\\
			$41$ & $21$ & 1.2 & 1.36 & 0.3 & 0.98 & 1.43 & 1.21 & $0.12$ & $38$ & $-5$\\
			$41$ & $21$ & \textbf{0.81} & 1.13 & 1.53 & 3.1 & 1.36 & 1.39 & $2.26$ & $30$ & $48$\\
			$41$ & $21$ & 3.49 & 3.47 & 9.35 & 6.0 & 6.92 & 3.81 & $1.47$ & $243$ & $-8$\\
			$42$ & $21$ & 2.46 & 2.76 & \textbf{2.29} & 4.12 & 3.68 & 2.65 & $5.64$ & $163$ & $21$\\
			$42$ & $21$ & 1.18 & 0.52 & 0.53 & 0.87 & \textbf{0.18} & 1.15 & $0.9$ & $17$ & $42$\\
			$42$ & $21$ & 5.89 & 5.86 & \textbf{2.3} & 3.05 & 6.19 & 4.12 & $3.71$ & $267$ & $5$\\
			$42$ & $21$ & 0.72 & 0.73 & 3.0 & 5.11 & \textbf{0.21} & 0.29 & $3.72$ & $16$ & $219$\\
			$42$ & $21$ & 2.76 & 2.2 & 4.95 & 3.75 & 2.68 & 3.05 & $0.56$ & $230$ & $-7$\\
			$42$ & $21$ & 3.31 & 3.31 & 7.17 & 5.41 & \textbf{3.13} & \textbf{3.13} & $5.81$ & $327$ & $8$\\
			$42$ & $21$ & \textbf{0.62} & 1.21 & 0.69 & 0.71 & 0.81 & 0.65 & $6.84$ & $62$ & $100$\\
			$42$ & $21$ & 3.46 & 3.5 & \textbf{1.21} & 3.95 & 1.66 & 1.38 & $4.01$ & $129$ & $22$\\
			$42$ & $21$ & 5.0 & 5.05 & 3.98 & 3.12 & 2.42 & 1.57 & $0.24$ & $101$ & $-13$\\
			$42$ & $21$ & 3.36 & 3.39 & 0.64 & \textbf{0.46} & 1.94 & 2.36 & $7.06$ & $65$ & $102$\\
			$43$ & $21$ & 0.62 & 0.61 & 4.35 & 1.89 & 2.45 & \textbf{0.26} & $3.84$ & $18$ & $199$\\
			$43$ & $21$ & 5.16 & 5.18 & 5.17 & 4.41 & \textbf{1.78} & 3.36 & $31.83$ & $161$ & $187$\\
			$43$ & $21$ & 3.06 & 3.03 & \textbf{0.13} & 1.26 & 0.39 & 0.41 & $3.46$ & $14$ & $238$\\
			$43$ & $21$ & \textbf{0.65} & \textbf{0.65} & 2.25 & 2.36 & 2.01 & 1.21 & $1.42$ & $16$ & $48$\\
			$43$ & $21$ & 0.97 & 1.75 & 0.86 & \textbf{0.18} & 1.31 & 0.82 & $9.4$ & $12$ & $768$\\
			\bottomrule
		\end{tabular}
		&\quad
		\begin{tabular}[t]{ccccccccccc}
		\toprule
			{n} & {m} & {mviol} & {aviol} & {pc4} & {pc8} & {freq4} & {freq8} & {gurobi} & {nq} & {qal}\\
			\midrule
			$43$ & $22$ & 4.24 & \textbf{2.65} & 10.32 & 12.87 & 8.33 & 7.75 & $4.87$ & $235$ & $9$\\
			$43$ & $22$ & 1.49 & 1.5 & 3.41 & 5.69 & \textbf{1.33} & 1.46 & $2.73$ & $67$ & $21$\\
			$43$ & $22$ & 1.31 & 1.32 & 3.19 & 4.27 & 1.39 & 1.48 & $1.01$ & $31$ & $-10$\\
			$43$ & $22$ & 2.98 & 2.76 & 7.59 & 16.77 & \textbf{1.14} & 3.02 & $9.04$ & $98$ & $81$\\
			$43$ & $22$ & 1.12 & 1.12 & 0.57 & 0.99 & \textbf{0.31} & 0.41 & $4.83$ & $23$ & $197$\\
			$44$ & $22$ & \textbf{1.12} & 1.13 & 2.96 & 2.66 & 2.31 & 2.07 & $8.21$ & $26$ & $273$\\
			$44$ & $22$ & \textbf{1.96} & 1.97 & 4.45 & 3.21 & 2.33 & 3.85 & $11.34$ & $44$ & $213$\\
			$44$ & $22$ & 1.11 & 1.12 & 1.36 & 2.09 & 2.43 & 2.64 & $0.44$ & $17$ & $-39$\\
			$44$ & $22$ & \textbf{2.16} & 7.15 & 4.79 & 4.39 & 6.04 & 7.54 & $12.29$ & $316$ & $32$\\
			$44$ & $22$ & 12.17 & 8.14 & 4.0 & \textbf{3.91} & 4.08 & 4.92 & $4.23$ & $483$ & $1$\\
			$44$ & $22$ & 6.85 & \textbf{5.1} & 13.98 & 6.5 & 9.0 & 7.57 & $44.5$ & $512$ & $77$\\
			$44$ & $22$ & 5.38 & 5.81 & 3.19 & 6.06 & 3.76 & \textbf{2.51} & $54.8$ & $212$ & $247$\\
			$44$ & $22$ & 23.14 & 12.78 & \textbf{9.75} & 9.96 & 16.87 & 16.65 & $33.91$ & $1758$ & $14$\\
			$44$ & $22$ & 1.75 & 3.48 & 0.7 & \textbf{0.59} & 0.8 & 3.04 & $6.93$ & $39$ & $163$\\
			$44$ & $22$ & 2.55 & 1.11 & 3.26 & 3.8 & 4.81 & 3.15 & $1.09$ & $16$ & $-1$\\
			$45$ & $22$ & 1.14 & 1.14 & 2.71 & \textbf{1.05} & 1.55 & 1.73 & $12.2$ & $72$ & $155$\\
			$45$ & $22$ & 17.7 & \textbf{12.29} & 16.21 & 14.56 & 17.54 & 12.77 & $39.59$ & $937$ & $29$\\
			$45$ & $22$ & 1.3 & 1.31 & 3.56 & 3.19 & \textbf{0.33} & 0.44 & $7.72$ & $31$ & $238$\\
			$45$ & $22$ & \textbf{0.77} & 1.36 & 3.18 & 4.96 & 1.82 & 1.77 & $9.24$ & $20$ & $424$\\
			$45$ & $22$ & 1.15 & 1.15 & 6.14 & 2.56 & 2.24 & \textbf{0.28} & $6.93$ & $20$ & $332$\\
			$45$ & $23$ & \textbf{1.12} & 1.14 & 2.4 & 1.42 & 1.85 & 3.7 & $37.25$ & $21$ & $1720$\\
			$45$ & $23$ & \textbf{3.35} & 7.57 & 9.57 & 5.95 & 6.95 & 4.47 & $6.07$ & $216$ & $13$\\
			$45$ & $23$ & \textbf{0.8} & 2.31 & 8.28 & 5.03 & 3.71 & 2.61 & $12.38$ & $81$ & $143$\\
			$45$ & $23$ & 8.23 & 16.88 & 7.36 & 7.6 & 7.09 & 6.37 & $1.02$ & $533$ & $-10$\\
			$45$ & $23$ & 11.71 & 13.73 & 13.83 & 23.29 & 18.24 & \textbf{9.3} & $10.72$ & $313$ & $5$\\
			$46$ & $23$ & 2.83 & 2.82 & 0.81 & 7.71 & \textbf{0.24} & 0.3 & $10.67$ & $24$ & $435$\\
			$46$ & $23$ & 6.13 & 6.13 & 10.27 & 19.55 & \textbf{3.84} & 7.77 & $4.89$ & $196$ & $5$\\
			$46$ & $23$ & 5.71 & 3.6 & 0.61 & \textbf{0.43} & 1.62 & 2.11 & $47.42$ & $116$ & $405$\\
			$46$ & $23$ & 9.36 & \textbf{2.89} & 10.26 & 9.09 & 14.61 & 8.6 & $6.57$ & $216$ & $17$\\
			$46$ & $23$ & 2.18 & \textbf{1.7} & 4.75 & 2.27 & 3.23 & 3.1 & $5.92$ & $119$ & $35$\\
			$46$ & $23$ & 11.44 & \textbf{5.14} & 17.04 & 14.56 & 15.62 & 13.09 & $13.48$ & $577$ & $14$\\
			$46$ & $23$ & \textbf{2.62} & 2.68 & 5.92 & 4.86 & 5.56 & 3.06 & $3.65$ & $50$ & $21$\\
			$46$ & $23$ & 2.48 & 2.47 & 4.0 & 6.3 & 2.34 & \textbf{1.41} & $18.59$ & $88$ & $195$\\
			$46$ & $23$ & 2.22 & 2.2 & \textbf{0.35} & 1.48 & 0.57 & 0.62 & $67.97$ & $82$ & $825$\\
			$46$ & $23$ & \textbf{1.93} & 1.94 & 6.88 & 7.05 & 4.57 & 2.69 & $28.22$ & $36$ & $730$\\
			$47$ & $23$ & 1.7 & \textbf{1.69} & 6.56 & 3.68 & 2.74 & 1.72 & $183.57$ & $49$ & $3712$\\
			$47$ & $23$ & \textbf{0.96} & 1.0 & 3.63 & 2.13 & 3.64 & 1.68 & $60.84$ & $42$ & $1426$\\
			$47$ & $23$ & \textbf{4.93} & 9.57 & 25.88 & 32.14 & 38.77 & 38.04 & $24.99$ & $127$ & $158$\\
			$47$ & $23$ & 1.26 & 1.27 & \textbf{0.1} & \textbf{0.1} & 0.23 & 0.29 & $36.9$ & $9$ & $4089$\\
			$47$ & $23$ & 10.82 & 6.06 & 2.09 & \textbf{0.21} & 9.57 & 6.22 & $15.56$ & $17$ & $903$\\
			$47$ & $24$ & 4.85 & 4.84 & 10.89 & 8.64 & 5.68 & \textbf{2.3} & $3.57$ & $57$ & $22$\\
			$47$ & $24$ & 5.84 & 5.86 & 3.09 & 5.85 & \textbf{1.94} & 4.39 & $14.91$ & $145$ & $89$\\
			$47$ & $24$ & \textbf{1.95} & \textbf{1.95} & 11.86 & 11.32 & 5.99 & 5.11 & $37.78$ & $59$ & $607$\\
			$47$ & $24$ & 17.5 & 14.05 & 39.38 & 18.33 & 18.66 & 22.22 & $5.19$ & $139$ & $-64$\\
			$47$ & $24$ & 2.88 & 2.91 & \textbf{2.85} & 8.17 & 11.95 & 4.41 & $17.78$ & $123$ & $121$\\
			$48$ & $24$ & 1.63 & 1.63 & 2.58 & 2.53 & 2.24 & 1.04 & $0.47$ & $39$ & $-15$\\
			$48$ & $24$ & 7.85 & 7.92 & 9.24 & 7.17 & 15.97 & \textbf{4.26} & $20.38$ & $174$ & $93$\\
			$48$ & $24$ & 10.22 & 10.48 & \textbf{6.99} & 7.55 & 17.15 & 9.63 & $28.7$ & $445$ & $49$\\
			$48$ & $24$ & 1.6 & 1.63 & 2.13 & 3.55 & 3.77 & \textbf{1.26} & $19.45$ & $43$ & $423$\\
			$48$ & $24$ & 13.47 & 8.18 & 13.45 & 14.14 & 11.45 & 8.37 & $6.46$ & $136$ & $-13$\\
			$48$ & $24$ & 11.0 & 11.93 & 5.21 & 6.9 & 12.4 & 7.69 & $4.39$ & $192$ & $-4$\\
			$48$ & $24$ & 3.0 & 3.1 & 2.84 & 4.48 & 4.07 & \textbf{1.36} & $241.22$ & $64$ & $3748$\\
			$48$ & $24$ & 3.52 & 3.53 & 2.76 & 2.46 & \textbf{1.37} & 2.45 & $256.81$ & $102$ & $2504$\\
			$48$ & $24$ & \textbf{3.69} & 3.75 & 11.3 & 12.38 & 10.97 & 5.15 & $36.82$ & $73$ & $454$\\
			$48$ & $24$ & 2.03 & 2.03 & 1.48 & 1.46 & \textbf{1.12} & 1.53 & $58.25$ & $112$ & $510$\\
			$49$ & $24$ & \textbf{9.72} & 11.32 & 28.36 & 23.65 & 10.09 & 15.96 & $25.18$ & $286$ & $54$\\
			$49$ & $24$ & 4.51 & 4.28 & 2.52 & 2.75 & 2.56 & \textbf{0.64} & $560.06$ & $25$ & $22377$\\
			$49$ & $24$ & 5.06 & \textbf{4.52} & 10.1 & 7.95 & 9.24 & 8.65 & $525.87$ & $116$ & $4494$\\
			$49$ & $24$ & 2.96 & 2.94 & 6.15 & 6.65 & 2.24 & \textbf{1.9} & $72.89$ & $62$ & $1145$\\
			$49$ & $24$ & \textbf{2.32} & \textbf{2.32} & 10.13 & 4.23 & 6.18 & 4.86 & $49.47$ & $44$ & $1072$\\
			$49$ & $25$ & 3.43 & 3.66 & 13.8 & 6.77 & 7.03 & \textbf{3.23} & $467.8$ & $127$ & $3658$\\
			$49$ & $25$ & 7.34 & 6.94 & \textbf{5.3} & 5.83 & 6.47 & 8.2 & $42.55$ & $577$ & $65$\\
			$49$ & $25$ & 9.61 & 5.13 & 20.6 & 8.41 & 3.98 & \textbf{3.69} & $94.45$ & $161$ & $564$\\
			$49$ & $25$ & 5.98 & \textbf{2.67} & 8.96 & 9.09 & 7.62 & 3.98 & $89.12$ & $95$ & $910$\\
			$49$ & $25$ & 2.79 & 2.79 & 8.64 & \textbf{0.21} & 0.86 & 1.56 & $228.04$ & $12$ & $18986$\\
			$50$ & $25$ & \textbf{4.34} & 5.64 & 7.95 & 7.04 & 5.65 & 6.67 & $284.31$ & $197$ & $1421$\\
			$50$ & $25$ & 2.11 & 2.12 & \textbf{0.23} & 7.31 & 4.07 & 3.22 & $454.95$ & $17$ & $26748$\\
			$50$ & $25$ & 14.14 & 14.09 & 31.87 & 4.3 & 4.35 & \textbf{3.38} & $211.75$ & $175$ & $1191$\\
			$50$ & $25$ & \textbf{6.49} & 7.43 & 9.14 & 10.86 & 10.31 & 10.24 & $373.16$ & $135$ & $2716$\\
			$50$ & $25$ & \textbf{3.76} & 4.34 & 17.29 & 26.28 & 30.06 & 10.99 & $69.83$ & $148$ & $446$\\
			$50$ & $25$ & \textbf{1.74} & \textbf{1.74} & 6.6 & 5.86 & 4.63 & 5.37 & $188.75$ & $42$ & $4453$\\
			$50$ & $25$ & \textbf{1.57} & 1.58 & 4.6 & 4.22 & 3.61 & 4.82 & $23.5$ & $21$ & $1044$\\
			$50$ & $25$ & 4.68 & 4.72 & \textbf{2.47} & 41.96 & 13.16 & 30.99 & $25.32$ & $271$ & $84$\\
			$50$ & $25$ & 4.15 & 6.59 & 1.91 & 1.99 & \textbf{1.28} & 1.31 & $41.05$ & $43$ & $925$\\
			$50$ & $25$ & 6.46 & \textbf{4.53} & 27.26 & 25.81 & 12.68 & 12.66 & $167.96$ & $237$ & $690$\\
		\bottomrule
		\end{tabular}
	\end{tabular}
\end{table}

In Tables \ref{tab:t1} and \ref{tab:t2}, we report results from computational
experiments performed on the group of test instances, evaluating each of the
different branching strategies. In these tables the columns {\tt mviol} and {\tt
aviol} correspond to Algorithms \ref{alg:scmvb} and \ref{alg:allmvb},
respectively. The columns {\tt pc4} and {\tt pc8} correspond respectively to the
pseudo-cost $4$-look-ahead and $8$-look-ahead strategies, and {\tt freq4} and
{\tt freq8} correspond respectively to the frequency-based $4$-look-ahead and $8
$-look-ahead strategies. Table \ref{tab:t1} gives the number of nodes that the
branch-and- bound algorithm requires when using each of the branching
strategies. The final column reports the number of nodes that the Gurobi
Optimizer used when solving the problem directly. In terms of the number of
nodes explored in the branch-and-bound tree, the most-violated-constraint and
all-violated-constraints satisfaction branching schemes (Algorithm
\ref{alg:scmvb} and Algorithm \ref{alg:allmvb}, respectively) are clear winners.

Table \ref{tab:t2} reports the time taken, in seconds, for each of the branching
strategies and the Gurobi Optimizer to solve the problem to optimality. The
number of queries to the UBQP oracle and the quantum annealing leniency (QAL) of
the quantum branch-and-bound algorithm with respect to the Gurobi Optimizer are
respectively given in columns {\tt nq} and {\tt qal}. The {\tt qal} column is
computed by taking the difference between the time taken by the best branching
strategy and the Gurobi Optimizer, and dividing by the number of queries. The
entries in this column can be viewed as the maximum threshold of the average
time, in milliseconds, to perform each quantum annealing process in order to
solve the original problem faster than the Gurobi Optimizer. Note that the
frequency-based branching heuristics terminate in the least amount of
computational time. A summary of the results of the comparison of these
different branching strategies is provided in Table
\ref{tab:branching-comparison}. In this table, the average run time is the
geometric mean of the columns of Table \ref{tab:t2}.

\setlength{\tabcolsep}{\defaulttabcolsep}
\renewcommand{\arraystretch}{1}
\begin{table}[hb]
	\caption{Comparison of branching strategies.}
	\label{tab:branching-comparison}
	\begin{tabular}{rl}
		\begin{tabular}[t]{r|cccccc}
		\toprule
		branching strategy & mviol & aviol & pc8 & pc8 & freq4 & freq8 \\
		\midrule
		average run time & 2.088 & 2.121 & 2.26 & 2.345 & 1.892 & 1.699 \\
		wins in node count & 75 & 78 & 5 & 8 & 8 & 19 \\
		wins in run time & 24 & 13 & 15 & 11 & 20 & 18 \\
		\bottomrule
		\end{tabular}
	\end{tabular}
\end{table}

\begin{wrapfigure}{t}{0.52\textwidth}
  \begin{center}
    \includegraphics[width=0.52\textwidth]{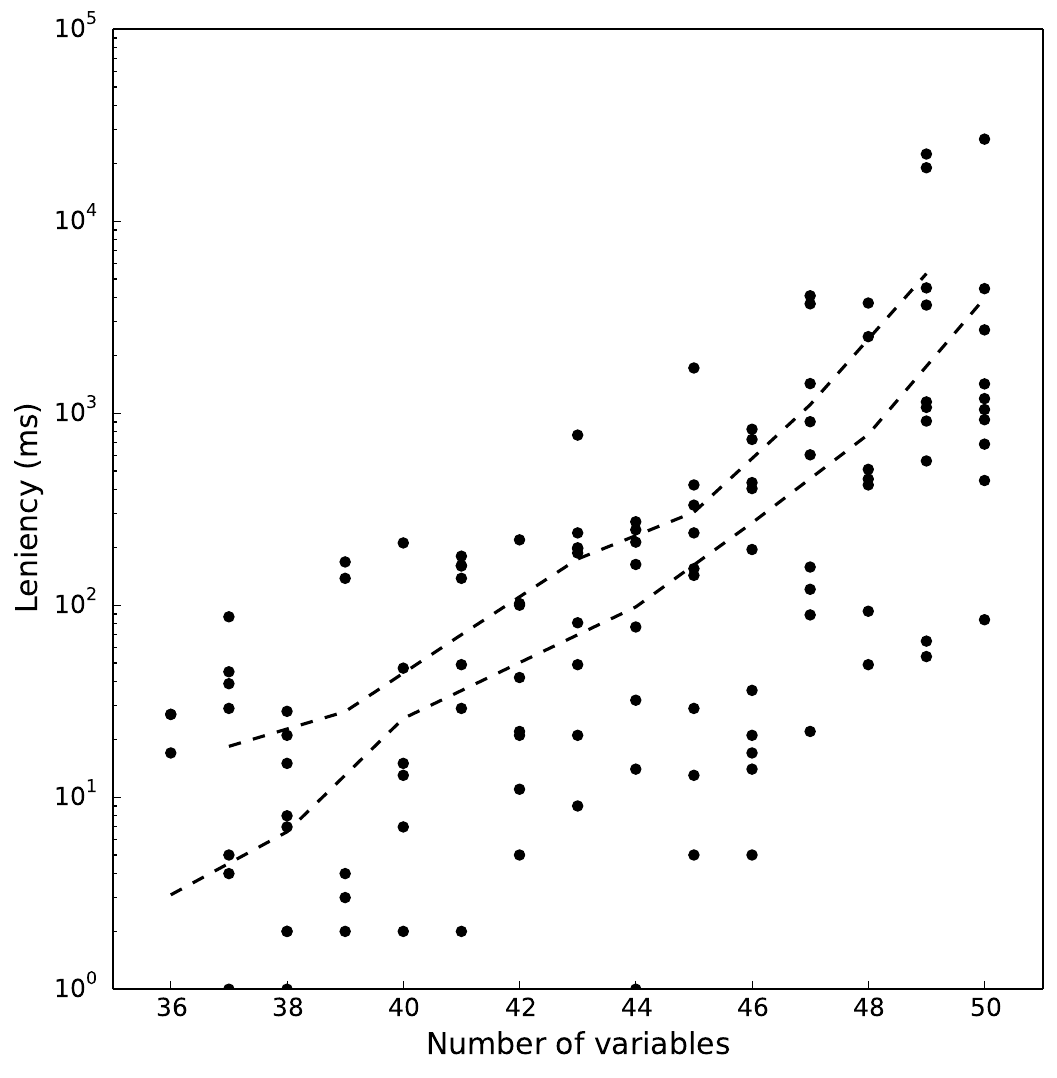}
  \end{center}
        \caption{\small{Quantum Annealing Leniency versus Problem Size}}
        \label{fig:2}
\end{wrapfigure}

In Figure \ref{fig:2}, we graph the QAL values from Table \ref{tab:t2} versus
problem size using a logarithmic scale. The dotted lines plot the mean values of
QAL. One dotted line corresponds to even values of $n$ and the other dotted line
corresponds to odd values of $n$.

This graph suggests an exponential growth in QAL with respect to the problem
size. We interpret this as an indication that the difference between the
computational time required for our CBQP approach (in conjuction with a scalable
UBQP oracle) and the Gurobi Optimizer grows exponentially with the size of the
problem.

\section{Discussion} 
\label{sec:Discussion} 

In this section we consider the specifics of the D-Wave devices as a physical
manifestation of a UBQP oracle. We imagine any implementation of quantum
adiabatic computing would have similar limitations, so we expect our algorithms
to be beneficial in overcoming them.

Quantum adiabatic devices have not thus far allowed for fully connected systems
of quantum bits. Due to this sparsity in the manufactured chips, the use of such
quantum computers requires solving a \emph{minor-embedding} problem, described
in the following section, prior to programming the chip according to appropriate
couplings and local fields \cite{AidanRoy2014}.

\subsection{Efficient embedding} 

Given a UBQP instance defined by a matrix $Q$, we define the underlying graph
$H$ as follows: for each variable $x_i$, we associate a vertex $v_i\in V(H)$;
and for each nonzero entry $q_{ij}\neq 0$ of $Q$ with $i \neq j$, we let $v_i$
and $v_j$ be adjacent in $H$. The minor-embedding problem is the problem of
finding a function $\phi:V(H)\rightarrow 2^{V(G)}$, where $G$ is the graph
defined by the quantum chip (that is, vertices correspond to the quantum bits
and edges correspond to the couplings between them), such that
\begin{enumerate}[(i)]
	\item for each $x\in V(H)$, the subgraph induced by $\phi(x)$ in $G$ is
	  connected;
	\item $\phi(x)$ and $\phi(y)$ are disjoint for all $x\neq y$ in $V(H)$; and,
	\item if $x$ and $y$ are adjacent in $H$, there is at least one edge between
	  $\phi(x)$ and $\phi(y)$ in $G$.
\end{enumerate}
Note that for any induced subgraph of $H$, a minor embedding can be found simply
by restricting $\phi$ to the vertices of the subgraph. In our CBQP approach,
the constraints contribute only to linear terms in the Lagrangian relaxation of
the problem, and hence the embedding of a UBQP problem at any node in the
branch-and-bound tree can be found from the parent node by restricting the
domain of $\phi$. That is, our method requires solving the minor-embedding
problem only once, at the root node of the branch-and-bound tree.

\subsection{Efficient programming of quantum chips} 

In every node of the branch-and-bound tree, all Lagrangian relaxations generated
have identical quadratic terms and only differ from each other in linear terms.
This suggests that if reprogramming the quantum chip can allow for fast updates
of previous setups, then the runtime of UBQP oracle queries can also be
minimized.

\subsection{Error analysis} 

Quantum bits currently have significant noise. For arbitrary choices of initial
and final Hamiltonians, the eigenvalues in the energy spectrum of the evolving
Hamiltonian of the system may experience gap closures. Furthermore, the
measurement process of the solutions of the quantum adiabatic evolution has a
stochastic nature. Each of these obstacles on its own indicates that the
solutions read from the quantum system are often very noisy and, even after
several repetitions of the process, there is no guarantee of optimality for the
corresponding UBQP problem. In order to make our method practical, with a proof
of optimality, it is necessary to develop a framework for error analysis for the
quantum annealer. Note that for our purposes the solution errors can only
propagate to final answers in the branch-and-bound tree if the proposed lower
bound $\ell$ obtained at a node is greater than the actual lower bound $\ell -
\epsilon$, and the best known upper bound $u$ satisfies $u < \ell$. If this
situation occurs, then the proposed method incorrectly prunes the subtree rooted
at this node. This motivates the study of a framework of error analysis that can
provide a measure of certainty on the optimality of solutions of the UBQP oracle
in the above sense.

\section{Extension to quadratically constrained problems} 
\label{sec:Generalization}

It is straightforward to extend the method proposed here to quadratically
constrained quadratic programming (QCQP) problems in binary variables. In fact,
the Lagrangian relaxations of QCQP problems are also UBQP problems. The minor-
embedding problem to be solved at the root node takes the underlying graph $H$
as follows: for each variable $x_i$, we associate a vertex $v_i\in V(H)$; and
for any pair of distinct indices $i \neq j$, we let $(v_i,v_j)\in E(H)$ if and
only if the term $x_i x_j$ appears with a nonzero coefficient in the objective
function or in any of the quadratic constraints.

Assuming future quantum annealing hardware will allow higher-degree interactions
between quantum bits, more-general polynomially constrained binary programming
problems could also be solved using a similar approach via Lagrangian duality.

\section{Conclusions}

Motivated by recent advancements in quantum computing technology, we have
provided a method to solve constrained binary programming problems using this
technology. Our method is a branch-and-bound algorithm where the lower bounds
are computed using Lagrangian duality and queries to an oracle that solves
unconstrained binary programming problems.

The conventional branching heuristics for integer programming problems rely on
fractional solutions of continuous relaxations of the subproblems at the nodes
of the branch-and-bound tree. Our lower-bounding methods are not based on
continuous relaxations of the binary variables. In particular, the optimal
solutions of the dual problems solved at the nodes of the branch-and-bound tree
are not fractional. We proposed several branching strategies that rely on
integer solutions and compared their performance both in time and number of
nodes visited, and to the conventional branching strategies that rely on
fractional solutions.

To compare the performance of our algorithm using quantum hardware against a
classical algorithm, we introduced the notion of quantum annealing leniency.
This is, roughly speaking, the average time a query to the UBQP oracle would
have in order to remain as fast as the benchmark algorithm for solving CBQP
problems. In our computational experiments this benchmark algorithm is that of
the Gurobi Optimizer.

Finally, although our focus was on quadratic objective functions and linear
inequality constraints, in Section \ref{sec:Generalization} we discussed how
this method can be generalized to higher-order binary polynomial programming
problems.

\section*{Acknowledgements}

The authors would like to thank Mustafa Elsheikh for his spectacular
implementation of the experimental software code, and Majid Dadashi and Nam Pham
from the software development team at 1QBit for their excellent support. We are
thankful to Marko Bucyk for editing a draft of this paper, and also to
Tamon~Stephen, Ramesh~Krishnamurti, and Robyn~Foerster for their useful
comments. This research was funded by 1QBit and partially by the Mitacs
Accelerate program.

\end{document}